\documentclass[12pt]{article}

\usepackage{amsthm,extmath}
\usepackage{latexsym}

\newcommand{\be}{\begin{equation}}
\newcommand{\ee}{\end{equation}}
\newcommand{\ba}{\begin{eqnarray}}
\newcommand{\ea}{\end{eqnarray}}

\theoremstyle{plain}
\newtheorem{tm}{Theorem}
\newtheorem{pr}[tm]{Proposition}
\newtheorem{co}[tm]{Corollary}
\newtheorem{lm}[tm]{Lemma}

\theoremstyle{remark}
\newtheorem{rk}{Remark}
\newtheorem{example}[rk]{Example}

\theoremstyle{definition}
\newtheorem{df}{Definition}

\begin{document}

\def\thesubsection{\arabic{subsection}}
\def\thesubsubsection{\thesubsection.\arabic{subsubsection}}

\title{\bf Differential equations of order 1 in differential fields of zero
characteristic$^{1}$ }

\author{\bf Jerzy Stry\l{}a}

\date{}

\maketitle

\footnotetext[1]{e-mail address: jstryla@wp.pl}


\begin{abstract}

I begin from a particular field of generalised Puiseux series and
investigate a class of nonlinear differential equations in the
field. It is appeared that the main part of differential equation
determines solvability and positions of resonance i.e. the
appearance of a free constant in solutions. Secondly, for any
singular differential equation in a differential field of the form
$Q(y)\dot y=P(y)$, $P,Q$ polynomials. Then 
the greatest Picard-Vessiot extension exists.
It is shown that the arising
set of solutions is obtained from algebraic equations labelled by
constants. Sufficient conditions for the PV extension being
extended liouvillian are delivered.
\end{abstract}


\subsection{Preliminaries}
     Differential   Galois   theory  gives  a  complete  view  of
algebraic extension structure of solutions for a given linear
differential equation \cite{Zo99, Ma94}. For proceeding the
research of general, nonlinear equation

\be   \label{pq}
\dot{y}=\frac{P(y)}{Q(y)},
\ee
where $P$ and $Q$, polynomials under a field $K$ have no common
roots in algebraic closure of $K$, firstly I elaborate a special
field of generalised Puiseux series in the section. The existence
problem was partially solved in \cite{Ar97,Ar00} by the method of
Newton polygons. For other approaches see for example \cite{BT98}.
Section \ref{P} refines the result. The introducing contour of
differential equation defines position of solutions beginning and
resonance. For a general differential field the existence of a
maximal Picard-Vessiot extension (the greatest one with the
isomorphism identification) via Eq. (\ref{pq}) is proved in section
\ref{PV}, theorem \ref{max} and corollary \ref{min}. In the place
of investigation of the group of symmetries of the extension
\cite{Zo99} for the case of linear equations, here, one focuses on
a rational functions associating with (\ref{pq}), (\ref{nicform}).
In the way a sufficient condition for realizing the extension by
finite number of operations: successive adding an exponent of
integral, an integral (liouvillian extension) or algebraic
extension appears as theorem \ref{end}.

I begin from a definition of (generalised) Puiseux field.
Everywhere in the paper a field $K$ is of zero characteristic.
Therefore, $\bbold{Q}\subset K$. I will further assume that an
ordered $\bbold{Q}$-linear space with order $(R_K,\le)$ is given,
i.e. $v_1>v_2 \Rightarrow v_1+v>v_2+v$ and $\bbold{Q}\ni q>0\wedge
R_K\ni v>0\Rightarrow qv>0$. An example is $R_K=\bbold{Q}$ or
$R_K=\bbold{R}\cap K$. The last case is equivalent to extra
property $\forall (r,\epsilon \in R_K)\exists (N\in \bbold{N})\
(r>0
\wedge \epsilon > 0)\Rightarrow n\epsilon > r$.

The completing of $R_K$ to $\bbold{Q}$-space with linear order
$\bar{R}_K$ by Peano construction enriches the space with suprema
and infima of nonempty and properly bounded subsets. $\bar{R}_K$ is
also a completion with respect to the norm $|\cdot|:R_K\mapsto R_K$
such that $|r|=r$ for $r\ge 0$ and $|r|=-r$ for $r<0$.
\begin{df}
Let $S\subset 2^{R_K}$ be the family of all well ordered sets of
$(R_K,\le)$. The Puiseux field $K_P[[x]]$ under the field $K$ is
the sum $\bigcup (K- \{ 0 \} ) ^s$ under all $s\in S$. The
operations are defined by the identification of the nonzero
elements of $K_P[[x]]$ with series of the form
$c_0x^{\mu_0}+c_1x^{\mu_1}+\ldots $, where the set of indices $\{\
\mu_i\}\in S$ and the set of coefficients $c_i\in  K -\{0\}$.
\end{df}
The multiplication is properly defined because of the property of
finite decomposability of elements of $s$, $s\in S$, with respect
to the additive operation in $R_K$. The set of non-decomposable
elements of $s$ is referred to as $B(s)=\{\mu\in s|\mu=\mu_1+\mu_2
\Rightarrow (\mu_1=\mu\ \vee \ \mu_2=\mu )\}$. If $s_1,s_2\in S$ I
write $s_1+s_2$ for $s_1\cup s_2$, $s_1\cdot s_2$ in the place of
$\{r_1+r_2|r_1\in s_1\
\wedge \ r_2 \in s_2 \}$ and $[s_i]$ for
$s_i+s_i\cdot s_i+\ldots $. If $\forall (x\in s_i)\ x\ge 0$ and
$R_K\subset \bbold{R}$ then $[s_i]\in S$. Extending the monoid
$(S,+)$ to the group $\bar{S} \equiv \{s_1-s_2|s_1\in S \wedge
s_2\in S\}$ one obtains a ring $(\bar{S},+,\cdot )$. Let $G(s)$
means the subset of $s$, $s\in S$, of elements without the
proceeding element in $s$. I put $N(s):=k$ for all $s$ in $S$ such
that $G^k(s)=GG\ldots G(s)$ is a finite set or the image of a
divergent to infinity sequence of elements of $s$ and $k\in
\{0\}\cup\bbold{N}\cup\{+\infty \}$ is the smallest such number and
$N(-s)=N(s)$.

The set of series $K_P^0[[x]] $ is a subring of $K_P[[x]]$, where
$K_P^0[[x]]$ contains by definition indices only from $N_{-1}(0)$.
If, additionally, $R_K\subset
\bbold{R}$ then $K_P^0[[x]]$ is a subfield.

Finding solutions of algebraic equations with coefficients in
$K_P[[x]]$ may be realized explicitly. In the way one may state
\begin{tm}      \label{ac}
Let $K$ be an algebraically closed field of zero characteristic.
Then $K_P[[x]]$ is algebraically closed.
\end{tm}
\begin{proof}
Any equation of degree one has a solution in $K_P[[x]]$. Let all
equations of degree less then $N-1$, $N\ge 1$, have a solution in
the field. Now, I take an algebraic equation,
\be   \label{al}
w(y):= \sum _{i=0}^N \alpha _iy^i=0,
\ee
$\alpha _i\in K_P[[x]]$, $\alpha_N\alpha_0 \ne 0$ and
$\alpha_i=\alpha_0(i)x^{\nu_0(i)}+\ldots $ for $\alpha_i \ne 0$.
Let $y=c_0x^{\mu_0}+c_1x^{\mu_1}+\ldots $ fulfils the equation. The
{\it contour function} $\frak{f}: R_K\longmapsto R_K$ is defined by
$\frak{f}(x)=\min \{\nu_0(i)+i x|i=0,\ldots,N \}$. The nonempty
finite set of breaking points $\{ x_b \}$ of $\frak{f}$ determines
all starting $\mu_0$ admitted by the equation by $\mu_0=x_b$. The
associated with $x_b$ polynomial $p_b(c_0)=\sum_{i\in
B}\alpha_0(i)c_0^i$, where $B\subset \{0,\ldots,N\}$ counts
realizations of the minimum of $\frak{f}$ in $x_b$, defines initial
$c_0\ne 0$ admitted for $\mu_0=x_b$. One may observe that the
number of available pairs $(\mu_0,c_0)\in R_K
\times K^{\ast}$
including multiplicities of $c_0$ is equal to $N$. Really, if a
polynomial has all roots in $K_P[[x]]$ then each $(\mu_0,c_0)$
including multiplicities begins a solution.

Let's assume that one possesses $y_0=c_0x^{\mu_0}+\ldots \in
K_P[[x]]$. Then $y=y_0+\bar{y}$ defines a new variable $\bar{y}$.
Eq. (\ref{al}) may be rewritten to
\ba     \label{rewr}
\sum_{i=0}^N \beta_i\bar{y}^i=0,  \\
\beta_i=\sum_{l=i}^N {{l}\choose{i}} \alpha_l y_0^{l-i}=
    \frac{1}{i!} \frac{d^i}{dy}w|_{y=y_0}.
\ea
I will refer to $y_0\in K_P[[x]]$ as a {\it partial solution} iff
for all $\Delta y_0=y_0-y_0^{\sigma }$, where $y_0^{\sigma }\ne
y_0$ is the restriction of $y_0$ to the indices less or less or
equal to $\sigma $, the term $c_{\sigma}x^{\mu_{\sigma }}$ from
$\Delta y_0:=c_{\sigma}x^{\mu_{\sigma } }+\ldots$ is an acceptable
initial term of a solution obtained from the contour of the
appropriate Eq. (\ref{rewr}). In (\ref{rewr}) still
$\beta_N=\alpha_N\ne 0$. If $\beta_0=0$ then $\bar{y}=0$ delivers a
complete solution.

Now, I will construct a solution in the following way. On each step
I use the contour function for the redefined equation to choose the
largest $\mu $, $\mu=x_b$, accompanied by appropriate $c\in
K^{\ast}$, $p_b(c)=0$. By the move I kill the lowest index element
$c_{\sigma }x^{\sigma }$ in $\beta_0$. Therefore, in the next step
$\sigma '$, the lowest index element of $\beta'_0$, will be greater
then $\sigma $. A set $T$ of triples $(\sigma ,c,
\mu)$ is arisen, where $\{\sigma \}$ appears to be a well ordered set by
the method of the construction. The minimal $\sigma $ is equal to
the lowest index of $\alpha_0$. I will show that in each step
labelled by $\sigma $ I always take the additional $cx^{\mu}$ such
that $\mu > g_{\sigma}\equiv \sup \{\bar{\mu}|
(\bar{\sigma},c,\bar{\mu})\in T\ \wedge\ \bar{\sigma}<\sigma \}$.
Firstly, I assume that
 the supremum $g_{\sigma} $ is the maximum. The case is reduced to the problem:
having a partial solution $c_0x^{\mu_0}$ of (\ref{rewr}) arising as
in the construction find a continuation $c_0x^{\mu_0}+c_1x^{\mu_1}$
according to the above rule. Let $c_0x^{\mu_0}$ cancels
$c_{\sigma}x^{\sigma}$ in $\beta_0$ and the new constituent to
remove in $\beta_0'=\sum^N_{i=0}\beta_i(c_0x^{\mu_0})^i$ is equal
to $c_{\sigma'}x^{\sigma'}$, where $\sigma' > \sigma $. The lowest
index expressions arising as a consequence of introducing
$c_1x^{\mu_1}$ has the forms
\be
c_1p^{(1)}(c_0)x^{M+(\mu_1-\mu_0)},\ldots ,c_1^k
p^{(k)}(c_0)x^{M+k(\mu_1-
\mu_0)},
\ee
where $p$ is the polynomial from the contour for $c_0x^{\mu_0}$,
$p^{(i)}(y)=\frac{d^i}{dy^i}p(y)$ is its $i^{th}$ derivative,
$k=\deg p$ and $M=\frak{f}_0(\mu_0)$ with the appropriate for the
first term contour function $\frak{f}_0$. At least one of them is
nonzero. If $\mu_1$ goes to $\mu_0^{+}$ then the indices of the
term goes to $\sigma$. Therefore, $\mu_1>\mu_0$ claimed to cancel
$c_{\sigma'}x^{\sigma'}$ exists. In turn, let's assume that a
partial solution $y_0$ exists such that $g_{\sigma'}$ is not the
maximum. Let's consider $y_0+c_1x^{\mu_1}$, where
$\mu_1=g_{\sigma'}$ and $c_1\in K^{\ast}$. Again, the lowest index
of $\sum^N_{i=0}\beta_i(y_0+c_1x^{\mu_1})^i$ containing
$(\mu_0,c_0)$ is equal to $\sigma_1\le \sigma <\sigma$, where $c_1$
is treated as a variable. The proper cancellation may be done.

\end{proof}
The rational functions $K(x)$ are a subset of $ K[[x]][x^{-1}]
\subset K [[ x^ {\frac{1}{s}} ]]$, where $K[[x^{\frac{1}{s}}]]:=
\bigcup_{s\in \bbold{N}} K[[x^{\frac{1}{s}}]][x^{-
\frac{1}{s}}]$ is the (classical) Puiseux series field. The field appears to be
algebraically closed as one of consequences of theorem \ref{ac}.
\begin{co}
$K_P^0[[x]]$, $K_{\bbold{Q}}^0[[x]]$ and $K[[x^{\frac{1}{s}}]]$ are
algebraically closed fields, where $R_K\subset \bbold{R}$ and
$K_{\bbold Q}^0[[x]]$ is obtained as a restriction of $K_P^0[[x]]$
to the rational indices.
\end{co}
\begin{proof}
Let $w(y)=y^N+\alpha_{N-1}y^{n-1}+\ldots +\alpha_0=0$ be an
algebraic equation with coefficients in $K_P^0[[x]]$, $N>0$. Also,
let algebraic equations of degree less then $N$ with coefficients
from $K_P^0[[x]]$ are solvable in $K_P^O[[x]]$. Let $y_s\in
K_P[[x]]$ be a solution of $w$ and $\underline{\mu}=\limsup
\mu_n\in \bar{R}_K$, where $\{\mu_n\}_{n=1}^{\infty}$ is an
increasing sequence of the set of indices of $y_s$, be the minimal
number with the property. Following the theorem, in the series
$w(y_s)$ the converging indices appears initially via $p^{(i)}
(c_0)c_{\mu_n}^ix^{M+i(\mu_n-\mu_0)}$, where $p^{(i)}(c_0)\ne 0$
and $i=1,\ldots ,N$ is the smallest such derivative. But by the
rescaling $y\rightarrow x^{- M}y$ I may assume that $M=0$. Now, the
expressions can not find counterparts under the assumption.
Therefore, $y_s\in K_P^0[[x]]$.

The corollary for the restriction of indices to $\bbold{Q}$ follows
from the construction of solutions in theorem \ref{ac}.

The final case is the coefficients in $K[[x^{\frac{1}{s}}]]$. At
the moment we know that a solution $y_s$ is at least from
$K_{\bbold{Q}}^0[[x]]$. I make the proof of the part by induction
with respect to degree of algebraic equation. Let
$\lim_{n\rightarrow\infty }\mu_n=
\infty$, where $\{ \mu_n \}$ are ordered by $\bbold{N}$ indices of $y_s$. Otherwise,
the set of the indices is finite and $y_s\in K[[x^{\frac{1}{s}}]]$.
Also, let $w^{(i)}(y_s)\ne 0$ for all $i=1,\ldots,N-1$. It implies
that the contours appearing in the construction of the solution for
$n>N_0$, $N_0\in \bbold{N}$, differ only by horizontal lines
defined by $\sigma_n$. Moreover, theorem \ref{ac} guarantees that
all partial solutions of the algebraic equations have a
continuation to a complete solution, so all of them are subsets of
a solution. The number of solutions is finite, so one finds that
for $n>N_1$ the solution is constructed via right choice of
breaking points as in the proof of the theorem and then the right
breaking points of the contours, $\mu_n$, $n>N_2$, are the
intersection of the horizontal lines and the line $r\rightarrow
\nu_0(1)+r$. Otherwise, an infinite number of different partial
solutions beginning different solutions exists. Let
$N=\max\{N_0,N_1,N_2\}$. All coefficients and the partial solution
$y_s|_N$ have indices in $\frac{1}{s}\bbold{N}$ for $s\in
\bbold{N}$. The procedure of extension behaves the pointed set
indices. The proof is done.
\end{proof}
\begin{rk}
Using the transformations $y\rightarrow x^M y$ or $y\rightarrow
y^{-1}$ one finds that the following situation is generic in the
space of leading indices $\{\nu_0(i)\}_{i=0}^N$ of coefficients.
Namely, $\alpha_N=1$, $\nu_0(N-k)<0$ for a $N-k\in\{1,\ldots , N
\}$ and $\nu_0(i)\ge 0$ for $i\ne N-k$. Then the explicit solution
may be proposed according to \cite{Ma97,Ma96}:
\be  \label{formal}
y_s=(-a_k)^{\frac{1}{k}}-\sum_{i=0}^{\infty} \left(
\sum_{p=1}^{i+1}
\frac{1}{p!k^p} \prod_{j=1}^{p-1} (i-jk)\sum_{{i_1+\ldots+i_p=i+1}\atop {i_l\ne k}}
a_{i_1}\cdot\ldots\cdot a_{i_p} \right) (-a_k)^{-\frac{i}{k}},
\ee
where $a_i=\alpha_{N-i}$, $R_K\subset \bbold{R}$ and I assume all
conventions from the paper. The formal series on finite intervals
reduce under the generic restrictions to an algebraic expressions
with the $k$-root. Then the theorem and the corollary are an
immediate consequence of (\ref{formal}).
\end{rk}
For a further context it is essential that the schedule of
subspaces
\be    \label{subsets}
K[[x^{\frac{1}{s}}]]\subset K_ {\bbold{Q}} ^0 [[x]]\subset
K_{\bbold{Q} } [[x]] \subset K_P[[x]] \supset K_P^0[[x]]
\ee
are defined only by referring to properties of indices and,
additionally, they are closed for $[\mu ]$ operation for
nonnegative indices $\mu $ for $R_K\subset \bbold{R}$.


Keeping the extra assumption and enriching the Puiseux field with
the standard differentiation
\be \label{defdiff}
(c_ix^{\mu_i})^{\bullet}:= c_i\mu_i x^{\mu_1 -1 }
\ee
I may return to Eq. (\ref{pq}). The $\frac{P(y)}{Q(y)}\in
K_P[[x]](y)$ I will interpret as an element of
$K_P[[x]]_{\bbold{Q}}[[y]]$. The change may leads to limitations
for admitted substitution of $y$. Nevertheless, the whole initial
domain $K_P[[x]]-\{ y| Q(y)=0 \}$ may be realized by considering
transformations $y\longmapsto y-y_0$, $Q(y_0)\ne 0$, to reach all
solutions of (\ref{pq}). In the next section I replace the problem
of solving the equation by a new one:
\be    \label{PQ}
\dot{y}=f_{ij}x^{\nu_i}y^{\sigma_j}
\ee
with the additional condition
\be  \label{condition}
f\in K_P[[x]]_{\bbold{Q}}[[y]]\cap K_{\bbold{Q}}[[y]]_P[[x]]
\ee
and $f_{i\cdot}\ne 0$, $f_{\cdot j}\ne 0$. I assume that a choice
of branches of $y^{\frac{p}{q}}$ is given.

\subsubsection*{Language of valuation}    \label{valuation}
The presented analysis and the approach of the next section
suggests a general view. Let $K\varsubsetneq L$ be an algebraically
closed subfield of a field $L$. Then a minimal valuation ring $O$
exist such that $K\subset O\varsubsetneq L$. Then the field $O/O_+$
may be identified with $K$, where $O_+$ is the maximal ideal of
$O$. Really, the maximal subfield $\tilde{K}$ in $O-O_+\supset K$
is isomorphic to $O/O_+$.
If $\tilde{K}\ne K$ then a valuation ring $\tilde{O}$ in
$\tilde{K}$ such that $K\subset \tilde{O}\varsubsetneq \tilde{K}$
exists and the ring generated by $\tilde{O}, O_+$ is a valuation
ring in $K$ smaller then $O$. The contradiction.

The valuation function $v:(L^{\ast},\cdot)\mapsto (R_{K},\oplus )$
is the groups homomorphism defined by $O$, where $(R_{K},\oplus
)\equiv (L/O^{\ast},\cdot)$ has no nilpotent elements and has the
linear order structure $\le$ such that $k_1\le k_2\iff k_1 O\supset
k_2O$ for $k_1,k_2\in R_K$. Additionally, $k_1\le k_2 \Rightarrow
k_1\oplus k\le k_2\oplus k$ for each $k\in R_K$. Such linearly
ordered group is naturally extendible treated as a
$\bbold{Z}$-module to $\bbold{Q}$ linear space with order
$(R_K,\oplus,\cdot,\le)$. For example, if $L$ is algebraically
closed then the multiplication by $q\in \bbold{Q}$ in $R_K$ is
already present and is defined by $k O^{\ast}\mapsto k^q O^{\ast}$.

Now, let $L$ be algebraically closed and $\{b_i>0\}$ be a basis of
$R_K$ and a choice $\{x(b_i)\}\subset L$ is given such that
$v(x(b_i))=b_i$. The $\{x(b_i)\}$ have no relations, therefore,
$K(\{b_i\})\subset L$. The referring to
$x(b_{i_1})^{q_1}\cdot\ldots \cdot x(b_{i_n})^{q_n}$ by $x^r$,
where $q_i\in \bbold{Q}$ and $r=q_1b_1+\ldots+q_nb_n\in R_K$, is
assumed. Let $K_P[[x]]$ be the Puiseux series with the set of index
equals to $R_K$. Then $L\subset K_P[[x]]$. For a proof, it is
enough to show that any nontrivial extension of $K_P[[x]]\supset
\bar{O}$ leads necessary to set theory extension of $R_K$, where
$\bar{O}$ is the valuation ring of $v(k)=\mu_0$.
\begin{pr}
Let $K$ be an algebraically closed field, $R_K$ a
$\bbold{Q}$-linear space with a linear order and
$K_P[[x]](y)\supset \frak{O}$ be an extension of $K_P[[x]]\supset
\bar{O}$. Additionally, let $K=\frak{O}/\frak{O}^{\ast}$. Then
$K(y)=K$ or $R_K\varsubsetneq R'_K\equiv
K_P[[x]](y)/\frak{O}^{\ast}$.
\end{pr}
\begin{proof}
If $v(y-k)=0$ for each $k\in K_P[[x]]$ then $K(y)\subset K$.
Therefore, let $v(y)>0$. The rest cases are equivalent. I assume
that $v(y)\in R_K$. Let $v(x^{\mu_0})=v(y)$ and $\mu_0\in R_K$. If
$v(y-c_0x^{\mu_0})>v(y)$ for $c_0\in K^{\ast}$ then $c_0$ is
unique. If the procedure has no any halt point then
$y=c_0x^{\mu_0}+c_1x^{\mu_1}+\ldots$ in $v$-topology (i.e. $\lim
y_n=y \iff \forall (r\in R_K) \exists (N\in \bbold{N}) \forall
(n>N)\ y-y_n\in O_r $, $O_r:=\{k\in L|v(K)\ge r \}$). So, for
$\bar{y}:=y-k$, $k\in K_P[[x]]$, $v(\bar{y}-cx^{\mu})=\mu>0 $ for
all $c\in K$. It implies that $v(y/x^{\mu}-c)=0$ and $y/x^{\mu}\in
\frak{O}/\frak{O}^{\ast}-\bar{O}/\bar{O}^{\ast}$.
A contradiction.
\end{proof}

If an injection of $\bbold{Q}$-linear space $i: R_K\mapsto K$
exists one may identify $R_K$ with a subset of $K$. In the case a
differentiation is defined by:
\be
(x^{\mu})^{\bullet}=\mu \frac{x^{\mu}}{x(b_1)},
\ee
where $b_1$ is an element of a chosen basis of $R_K$. An
identification $\bbold{Q}\cdot b_1$ with $\bbold{Q}$ via $b_1=1$,
assumption that $\dot{K}=\{0\}$ and the rule of independent
differentiation of each "monomial" in $k\in K_P[[x]]$ return us to
formula (\ref{defdiff}).

The inequality $v(\dot{y})-v(y)\ge 1$ becomes an equality for
$v(y)\ne 0$, where $y\in K_{\bbold{Q}}^0[[x]]$. In the approach the
case $v(y)=0$ is distinguished.

\subsection{Existing theorem and resonance} \label{P}
\subsubsection{Domain of differential equation}
Let's start from a characterization of the domain of (\ref{PQ}).
All $y=c_0x^{\mu_0}+\ldots $, where $\mu_0>0$, are properly
substituted to the equation, compare proposition \ref{D0} below.
Nonzero constants from the domain constitute the set $D_0\equiv
\{c\in K|\forall(i\in \nu)\ f_i(c)<\infty
\}$. It is reasonable to consider also a class of fields $K^{ns}$
being also a normalized space over $\bbold{Q}$, for example
$K^{ns}=\bbold{C}$. In the way, the condition $c\in D_0$ reads as
\be  \label{D_0}
\forall_i \sum_j f_{ij}c^{\sigma_j}< \infty
\ee
admits infinite summations. One defines $D_0'\subset D_0$ by
\be                 \label{D_0'}
c_0\in D_0' \iff \forall_{i\in \nu}\ \limsup_{j\rightarrow \infty }
\sqrt[j]{|f_{ij}c_0^{\sigma_j}|}<1 \}.
\ee
For $y=c_0+\ldots $, $c_0 \in K^{ns}-\{0\}$, the condition $c_0\in
D_0$ is necessary for staying in the domain. It is appeared that
$c_0\in D_0'$ is sufficient. We have
\begin{pr}     \label{D0}
Let $D_f\equiv \{y\in K_P^{ns}[[x]]|f(y)\in K_P^{ns}[[x]] \}$ be
the domain of the differential equation. Then all $y=c_0+\ldots$,
$c_0\in D_0'$ and $y=c_0x^{\mu_0}+\ldots $, $\mu_0>0$, belong to
$D_f$.
\end{pr}
\begin{proof}
We start from the observation that
\be
(c_0x^{\mu_0}+c_1x^{\mu_1}+
\ldots)^{\sigma_j}=c_0^{\sigma_j}x^{\sigma_j \mu_0}(1+
d_1(j) x^{\mu_1-\mu_0}+\ldots )
\ee
Therefore, the indices after substitution have the form
$\nu_i+\sigma_j
\mu_0+[\mu_k-\mu_0]$. Then for $\mu_0>0$ the sums of coefficients before an
admitted index in $f(c_0x^{\mu_0}+c_1x^{\mu_1}+\ldots )$ are always
finite.
Let $\mu_0=0$. Then the convergence of the lowest term coefficient
in $f(y)$ is guaranteed by (\ref{D_0'}). The coefficients $d_k(j)$
of $x^{\delta\mu(k)}$, for $[\mu_i-\mu_0] \ni \delta\mu(k)>0$, are
polynomials of $c_{k(1)}/c_0,\ldots ,c_{k(l)}/c_0$ spanned by a
finite number of monomials $m_{n}(k):= \prod_i ( c_{k(i)}/c_0
)^{n_i}$ such that $\delta\mu(k)=\sum n_i \delta\mu(k(i))$, where
$i=1,2,\ldots, I_k $, $n_i, I_k\in \bbold{N}$ and
$\delta\mu(k(i))\in [\mu(i)] -\{ 0 \}$ are elements of admitted
decompositions of $\delta\mu(k)$. For $\delta\mu(l) \in
[\mu_i-\mu_0]-\{\mu_i-\mu_0 \}$ I put $c_{\delta\mu(l)}=0$. The
coefficients $\gamma_{n}(k,j)\in \bbold{Q}$ defined by
\be  \label{gamma}
d_k(j)=\sum \gamma_n (k,j) m_n (k)
\ee
may be calculated directly. For $\sigma_j\in \bbold{N}$
\be   \label{sigmanat}
\gamma_n (k,j)=\frac{\sigma_j !}{n_0!n_1!\ldots n_{\sigma_j}!},
\ee
where one introduces $n_0:=\max\{
\sigma_j - \sum_{l=1} n_l,1 \}$. For  $\sigma_j=-1$ the coefficient is a
number. For $\sigma_j=\frac{1}{q}$, $q\in \bbold{N}$ one need to
solve a finite number of algebraic equations defined by relations
arising among $x^{\delta\mu(k(i))}$. The well-defined number is a
rational function of $\sigma_j$ depending on expressions as in
(\ref{sigmanat}). In the same way coefficient arise for general
$\sigma_j\in\bbold{Q}$. Now, the summability of $\sum
f_{ij}c_0^{\sigma_j}d_k(j)$ keeping constant $\nu_i+\delta \mu_k$
is under consideration. Again (\ref{D_0'}) is sufficient for
convergence, because the only infinite sum may appear is through
$j$.
\end{proof}
\begin{rk}  \label{minus}
For considering any $K$ and any initial index one must assume that
the set $\{ \sigma_j \}$ is finite. Then $D_0=K$ and
$D_f=K_P[[x]]$.

For any $K$ and $\mu_0\ge 0$ admitted as the beginning index it is
enough that $\forall (i) \sum_jf_{ij}< \infty$ i.e. there are only
a finite number of $f_{ij}$ for each fixed $i$. In the situation
$D_0=K$ and $D_f$ contains all Puiseux series with nonnegative
indices.
\end{rk}

\subsubsection{Necessary conditions}
For investigating of necessary conditions for existence of
solutions of Eq. (\ref{PQ}) one defines the contour function for
the equation
\be
\frak{f}(x):= \min\{\nu_i+\sigma_j x, x-1|x\in R_K\wedge
f_{ij}\ne 0 \}.
\ee
The function's graph for $x\ge 0$ (or for $x\in R_K$, while
$\sigma$ is finite) is a broken line with finite number of breaking
points. Also for each $x> 0$ (resp. $x\in R_K$) a finite set of
pairs $(\nu_i,\sigma_j)$ exists such that
$\frak{f}(x)=\nu_i+\sigma_j x$. The following proposition is an
answer for the first step of the equation solving:
\begin{pr}  \label{necessary}
If $y\in K_P[[x]]$, $y\ne 0$, is a solution of (\ref{PQ}) then the
set of indices of $y$, $i(y)$, begins from $\mu_0$ such that
\begin{enumerate}
\item[(a)] $\mu_0=0$ or
\item[(b)] $\mu_0 $ is a breaking point of $\frak{f} $ and $\mu_0\ne 0$ or
\item[(c)] $\frak{f}(\mu_0)=\mu_0-1$
and $\mu_0=f_{ij}\ne 0$, $\nu_i=-1$, $\sigma_j=1$.
\end{enumerate}
\end{pr}
\begin{proof}
Let $\mu_0\ne 0$ and $y=c_0x^{\mu_0}+\ldots $ fulfils the
differential equation. Then the algebraic equation for $c_0$,
\be      \label{vertex}
p_{\mu_0}(c_0):= \sum_{i,j,\nu_i+\sigma_j \mu_0=\frak{f}(\mu_0)}
f_{ij}c_0^{\sigma_j}-\mu_0c_0=0,
\ee
admits the nonzero roots only for the cases (b) or (c).
\end{proof}
\begin{rk}
The choice of branches of $y^\sigma_j $ may sometimes cause Eq.
(\ref{vertex}) has no nonzero solutions. Further, taking $c_0\ne 0$
I will suppose a coherent initial definitions.
\end{rk}
The above proposition state that in the accepted area of initial
indices there are only finite number of possibilities. If there are
no breaking points of $\frak{f}$ then $y=0$ is a solution.
Moreover, it is appeared that $c_0$ is not defined in item (c) and
no (b) of the proposition. The same is true in (a) under
\be \label{Aa}
\nu_0+1>0
\ee
i.e. $c_0\in D_0$. If (\ref{Aa}) is not the case then the following
condition for $c_0$ immediately arises $\sum_j f_{0j
}c_0^{\sigma_j}=0$.

The distinction is close to a general classification of solutions.
The {\it algebraic type} solution of Eq. (\ref{PQ}) is such which
in first step of solving does not bind left hand side of the
equation. The rest of solutions are called ${\it proper}$ ones. In
the way (c) leads only to proper solutions, (a) without (\ref{Aa})
only to algebraic type and (b) refers to proper ones iff the
breaking point lies on $y=x-1$. From the properties of $\frak{f}$
it implies that the maximal initial indices' number of proper
solutions is $4$, in case (b) at most two.

It will appear that the free parameters $c_i$ may arise only one
time in a solution. For each $(\mu_0,c_0)$ from cases (b)(c) of
proposition \ref{necessary} I define the {\it resonant index }
$\mu_r$ by
\be \label{resonant}
\mu_r=\sum f_{ij}c_0^{\sigma_j-1}\sigma_j ,
\ee
where the summation is taken through $\{(i,j)|\nu_i+\sigma_j
\mu_0=\frak{f}(\mu_0)\}$. One finds that the case (c) begins from the resonant value
$\mu_0=\mu_r$. For terminological consistency case (a) restricted
to (\ref{Aa}) will also be called as resonant.

In turn, I will find a well ordered set covering indices of proper
solutions $y=c_0 x^{\mu_0}+ \ldots $.
\begin{pr}   \label{indices}
Let $R_K\subset \bbold{R}$ and $c_0x^{\mu_0}$, $\mu_0\ge 0$ and
$\neg (\mu_r > \mu_0)$, may start a proper solution. If
$y=c_0x^{\mu_0}+\ldots $ is a solution then the indices set $i(y)$
is a subset of $i_{\max}(c_0x^{\mu_0}):= [ \{ \nu+1 \} \cdot \{
\mu_0 (\sigma -1) \} ] \cdot \{ \mu_0 \}$.
\end{pr}
\begin{proof}
Let $y=c_0x^{\mu_0}+c_1x^{\mu_1}+\ldots $ be a solution of
(\ref{PQ}) written in the convention such that $\{\mu \}=[\mu]$ and
$c_i=0$ are admitted for $i\ne 0$. One defines a function
$\frak{f}_0(x) := \min \{\nu_i+\mu_j x|x\in R_K
\wedge f_{ij}\ne 0 \}$.
Let $\mu_0\ne 0$. Then
\be     \label{0}
\mu_0-1=\frak{f}_0(\mu_0)=\nu_i+(\sigma_j-1) \mu_0 +\mu_0
\ee
for some $\nu_i,\sigma_j$ realizing the minimum. Let $\delta
\mu_i:=
\mu_i-\mu_0$. Then $y^{\sigma_j}=c_0^{\sigma_j} x^{\mu_0 \sigma_0} (1+
\ldots+
d_i (j) x^{\delta \mu_i}+\ldots )$. A comparison of both sides of
(\ref{PQ}) pitches the condition for the indices:
\be      \label{eqonind}
\delta\mu_l= \nu_i+1  +(\sigma_j-1)\mu_0+  \delta\mu_k
\ee
with $\nu_i+1+ (\sigma_j-1) \mu_0\ge 0$ for all $i,j$ such that
$f_{ij}\ne 0$. The accompanying coefficient equation is the
following
\be  \label{eqoncoe}
c_l \mu_l = \sum f_{ij}c_0^{\sigma_j}d_j(k)=
\sum_{i,j, \frak{f} (\mu_0) = \nu_i+\sigma_j \mu_0  }
f_{ij} c_0 ^ {\sigma_j-1} \sigma_jc_l+ \sum (k<l).
\ee
Now, if in Eq. (\ref{eqonind}) $\delta\mu_l$ has no realization for
$k<l$ then (\ref{eqoncoe}) implies $c_l \mu_l=c_l \mu_r+0 $, so
$c_l=0$. Therefore, starting from (\ref{0}) and following by
(\ref{eqonind}) for $k\ne l $ one gains the thesis. The proof for $
\mu_0=0 $ follows similarly.
\end{proof}

\subsubsection{Existence theorems}
Now, I pass to the proper solutions' continuation problem.
\begin{tm}  \label{proper}
Let $R_K\subset \bbold{R}$ and $c_0x^{\mu_0}$, $\mu_0\ge 0$, is
admitted as an initial term of a proper solution. Then the
following statements are fulfilled:
\begin{enumerate}
\item[(i)]  if $\mu_0\ne 0$, $\neg (\mu_r > \mu_0)$ then a unique
continuation exists,
\item[(ii)]  if $\mu_0\ne 0$, $\mu_r>\mu_0$ and
$\mu_r \notin i_{\max}(c_0x^{\mu_0})$ then any $c_r\in K$
determines a continuation, where $c_r$ is the coefficient before
$x^{\mu_r}$,
\item[(iii)] if $\mu_r>\mu_0>0$ and $\mu_r\in i_{\max}(c_0x^{\mu_0})$  then
a finite number of steps solves the alternative: the continuation
is as in case (ii) or the solution is terminated (negative
resonance),
\item[(iv)] if $\mu_0=0$ then there is exactly one continuation .
\end{enumerate}
\end{tm}
\begin{proof}
Item (i) follows directly from the proof of proposition
\ref{indices}. Each $c_l\in K$ for $i_{\max}(c_0x^{\mu_0})\ni
\mu_l>\mu_0$ is uniquely countable from (\ref{eqoncoe}). In case
(ii) $\delta\mu_r:=\mu_r -\mu_0$ can not be of the form
(\ref{eqonind}). This is the reason why appropriate (\ref{eqoncoe})
becomes a tautology. Any $c_r\in K$ is admitted and the successive
solving may be continued from $\mu_r$ with additional
$c_rx^{\mu_r}$. The new set of admitted indices is equal to
$i_{\max}(c_0x^{\mu_0})+i_{\max}(c_0x^{\mu_0})\cdot [\mu_r-\mu_0]$.

The situation (iii) delivers an extra condition on the level of
$\mu_r$. The index $\mu_r$ defined a finite set of non-decomposable
elements $B_r\subset B(i_{\max}(c_0x^{\mu_0}))$ which may appear in
decomposition of $\mu_r$. According to (\ref{eqonind}) and
(\ref{eqoncoe}) the coefficient $c_{\mu}$, $\mu\in [B_r]$ and $\mu
< \mu_r$ are uniquely determined by a finite number of algebraic
equations. They meet themselves on the level $\mu_r$ in appropriate
(\ref{eqoncoe}) and the relation, non-containing $c_r$, needs to be
justify with respect to its logical value. For example, if
$\mu_r\in B(i_{\max}(c_0x^{\mu_0}))$ then one meets case (ii).
Finally, in the last case $\mu_0=0$ the lowest Eq. (\ref{eqonind})
has the form
\be  \label{Aa'}
\mu_1=\nu_0+1>\mu_0,
\ee
which implies that all coefficient equations have no any obstacles
in solving and $c_0\in D_0'$ determines uniquely the rest
coefficients.
\end{proof}
\begin{rk}
For finite $\sigma$ the beginning coefficients $\mu_0<0$ may be
translated to nonnegative one via redefinition $y\rightarrow x^My$,
$M\in R_K$, and then theorem \ref{proper} covers also the
situation.
\end{rk}
\begin{rk}
The theorem for a general $R_K$ may be realized similarly with
omitting the predefining of indices. Then the case $\neg(\mu_r >
\mu_0)$ and $\mu_r\notin R_K$ may be treated as a resonant case
after a extension of $R_K\subset R_K'\ni \mu_r$ such that $\mu_r >
\mu_0$.
\end{rk}
With help of the distinction (\ref{subsets}) a classification of
Eq. (\ref{PQ}) may be imposed. From proposition \ref{indices} and
the above theorem one concludes the following fact.
\begin{co}
Let $R_K\subset \bbold{R}$. If the sets of indices $\{ \nu_i \}$
and $\{ \sigma_j \}$ belong to a subfield of (\ref{subsets}) then
for the cases (i), (ii) with $c_r=0$, (iii) and (iv) the proper
solutions belong to the subfield as well.
\end{co}
It needs to be remarked that the proper solution with $\mu_0=0$
implies (\ref{Aa'}), so (\ref{Aa}), but not inversely. For all
choices of $c_0\ne 0$ such that $\sum_j f_{0j}c_0^{\sigma_j}=0$ we
are in algebraic type solutions. Nevertheless, for (\ref{Aa}) being
fulfilled theorem \ref{proper} is applied (item (iv)). A different
treating of the other cases is needed.
\begin{tm}   \label{algebraic}
Let $R_K\subset \bbold{R}$ and $c_0x^{\mu_0}$ is admitted as an
algebraic type solution and $\{\sigma
\}$ is finite. If $\mu_0=0$ I additionally assume that $\nu_0+1<0$ and
that $\frak{f}_0$ has a breaking point in $0$. Then the set of
continuations as a solution $C(c_0,\mu_0)$ in nonempty and
$\#C(c_0,
\mu_0)\le s2^{(\sigma_{\max}-\sigma_{\min})s-1}$, where $\sigma_{\max}:
=\max\{\sigma \}$, $\sigma_{\min}:=\min\{\sigma \}$ and $s\in \bbold{N}$
is the smallest common divisor of $\{\sigma \}$.
\end{tm}
\begin{proof}
Eq. (\ref{PQ}) has the form $\dot{y}=\sum_j p_j(x)
(y^{\frac{1}{s}})^{s \sigma_j}$, where $p_j(x):= \sum_i
f_{ij}x^{\nu_i}$, a branch of $y^{\frac{1}{s}}$ is fixed and $n_j:=
s\sigma_j \in \bbold{Z}$. Without lost of generality of the
consideration one may put $s=1$. The algorithm of solving begins
from the term such that $\frak{f}(\mu_0)<\mu_0 -1$ for $\mu_0\ne 0$
or $\mu_0=0$ is a breaking point of $\frak{f}_0$ with $\sum_j
f_{0j} c_0^{n_j}=0$. The next step is to solve the algebraic
equation
\be     \label{algebra0}
\sum_j p_j(x)y^{n_j}=0
\ee
for given $(\mu_0,c_0)$, see theorem \ref{ac}. Let $\mu_0\ne 0$
(the proof is analogous for $\mu_0=0$). The solution of
(\ref{algebra0}), $y_0$, will be up to index $\mu $, $\mu
-\mu_0=\mu_0-1 -\frak{f}(\mu_0)>0$, the solution of the
differential equation. Subsequently, I modify (\ref{algebra0}) to
\be \label{algebra1}
\sum_j p_j(x)y^{n_j}-\dot{y_0}=0.
\ee
Eq. (\ref{algebra1}) has a solution, $y_0'$, being a continuation
of $y_0|_{<\mu}$, see the proof of theorem \ref{ac}. Now, the
coincidence with the differential equation is improved up to $\mu'$
such that $\mu'-\mu_0=\mu -1 -\frak{f}(\mu_0)=2(\mu - \mu_0)$. The
cut series $y_0|_{<\mu},y_0'|_{<\mu'},\ldots $ become a sequence of
approximation of a solution of the differential equation. The
induction follows.

Each solution in the form $y^{\frac{1}{s}}$ may bifurcate at most
$(\sigma_{\max}-\sigma_{\min}-1)$ times. The numerical coefficient
is estimated.
\end{proof}
The construction offers the closeness of $K_{\bbold{Q}}^0[[x]]$,
$K_P^0 [[x]]$, $K_{\bbold{Q}}[[x]]$ for solving in algebraic type.
\begin{rk}
The first omitted case in the theorems about continuation is the
following: $\mu_0=0$, $\nu_0+1\le 0$ and $\frak{f}_0$ has no
breaking point in $0$. They have no continuation. For example,
equations having no solutions in $K_P[[x]]$ at all are of the form:
\be
\dot{y}= y^{\sigma}(f_0x^{-1}+\ldots),
\ee
where $\bbold{Q}\ni \sigma<0$ and $f_0\ne 0$. To solve them $\ln x$
function is needed.

The second and last case is $\mu_0=0$, $\nu_0+1=0$ and $\frak{f}_0$
has breaking point in $0$. Let $\{\sigma\}$ is finite. Then one may
redefine the differential equation by $y\rightarrow
x^{-\epsilon}y$, where $\epsilon >0$. The admitted, finite set of
initial terms $c_0$ appears as $c_0x^{\epsilon}$ in the new
equation, but, now, as a beginning of a proper solution (cases
(i)-(iii) from theorem \ref{proper}). Elementary contour analysis
shows that for each $c_0$ the resonance is placed in
$\mu_r=\sum_jf_{0j}c_0^{\sigma_j-1}\sigma_j=0$.
\end{rk}

\subsection{Picard-Vessiot extension} \label{PV}
An algebraic approach to Eq. (\ref{pq}) in a differential field (of
zero characteristic) $(K,
\dot{\ } )$ will be concentrated on existence of maximal set of solutions,
its algebraic structure and a decomposition of the extension into
liouvillian steps. It is assumed that the field of constants
$C\subset K$ is algebraically closed.

I begin from an observation that the algebraic closure
$\overline{K}\supset K$ is a differential field extension with the
uniquely defined extension of the differentiation. It behaves the
field of constants. Further, all extensions $(\widetilde{K},\dot{\
}) \supset (K,\dot{\ })$ such that the field of constants is kept
unchanged and that they are generated by solutions of (\ref{pq})
will be called Picard-Vessiot (PV) extensions. Firstly, let $K$ be
algebraically closed.
\begin{lm}  \label{kac}
Let $\{y_{sol}\}$ be the set of all solutions of (\ref{pq}) in
$K=\overline{K}$. Then the following statements are equivalent:
\begin{enumerate}
\item[(i)] the differential ring
\be
R:=
\left(
K[y]\left[ \frac{1}{Q} \right] \left[ \left\{
\frac{1}{y-y_{sol}}\right\}
\right], \dot{y}=\frac{P(y)}{Q(y)}
\right)
\ee
does not possess new invertible constants (NIC),
\item[(ii)] a nontrivial PV extension of $K$ exists.
\end{enumerate}
\end{lm}
\begin{proof}
(i) is a necessary condition of (ii). Therefore, it is enough to
explain the implication (i) $\Rightarrow $ (ii). The only nonzero,
prime, proper differential ideal of $R$ is of the form $(y-a)$,
$a\in K$. Then $R/I=K$, so $y$ is a solution from $K$. The
contradiction imply that $\{ 0 \}$ is the maximal differential
ideal. In the way $( K(y), \dot{y}=\frac{P(y)}{Q(y)} )$ is a PV
extension.
\end{proof}
\begin{pr}   \label{atleast}
Let $K=\overline{K}$ does not contain any solution of (\ref{pq}).
Then item (i) of lemma \ref{kac} is fulfilled.
\end{pr}
\begin{proof}
An invertible element of $R$ has the form
\be        \label{nicform}
\alpha (y-y_{l_1})^{k_1} \ldots (y-y_{l_n})^{k_n},
\ee
where $\alpha \in K^{\ast}$ and $y_{l_i}\in
\{y_{sol}\}\subset K
\vee Q(y_{l_i})=0$, and $k_i\in \bbold{Z}-\{ 0 \}$, and $i\ne j
\Rightarrow y_i\ne y_j$. It is a constant iff
\be  \label{nic}
\prod_j (y-y_{l_j}) \left(
\frac{\dot{\alpha}}{\alpha }Q(y)+\sum_{i} k_i\frac{-\dot{y_{l_i}}Q(y)+P(y)
}{y-y_{l_i}} \right) =0 .
\ee
I substitute $y=y_{l_i}$ to (\ref{nic}). One obtains that $y_{l_i}$ is a
solution.
\end{proof}
Let's assume the following notations. For a given differential filed $K$ $I_{\infty}\subset K[\{w_{1,i},w_{2,i},\dot{w}_{1,i},\dot{w}_{2,i},\ldots \}_{i=0}^{\infty}]\equiv R$ is the differential ideal in the differential ring $R$ generated by the equations:
\be  \label{formax}
e_N\equiv \sum_{i+j+k=N}\dot{w}_{1,i}w_{2,j}q_k-w_{1,i}\dot{w}_{2,j}q_k+
(i+1)w_{1,i+1}w_{1,j}p_k-w_{1,i}(j+1)w_{2,j+1}p_k=0.
\ee
Now, the proper statement has the following shape.
\begin{tm}  \label{max}
For any differential field $(K,\dot{\ })$ and Eq. (\ref{pq}) in the field 
the greatest PV extension exists.
\end{tm}
\begin{proof}
A maximal differential ideal $I_{\max}\supset I_{\infty}|_N$ of a 
differential ring
\be
K[1/(w_{1,k}w_{2,l}-w_{1,l}w_{2,k})] 
[\{w_{1,i},w_{2,i},\dot{w}_{1,i},\dot{w}_{2,i},\ldots \}_{i=0}^{N}]
\varsupsetneq I_{\infty}|_N
\ee
exists such that $(w_{1,k}w_{2,l}-w_{1,l}w_{2,k})^n\notin I_{\infty}|_N$ 
for $n\in {\bbold N}$ and $N\ge k,N\ge l$ and $|_N$ means a restriction
by $w^{(l)}_{t,i}=0$ for $t=1,2$, $i>N$ and $l\in {\bbold N}\cup \{0\}$. 
Really, let's assume that 
$(w_{1,k}w_{2,l}-w_{1,l}w_{2,k})^n\in I_{\infty}|_N$, $k\ne l$. Then 
$ (w_{1,k}w_{2,l}-w_{1,l}w_{2,k})^n=\sum e_i|_N r_i$, where $r_i\in R$. 
It implies that $ (w_{1,k}w_{2,l}-w_{1,l}w_{2,k})^n=
((k-1)w_{1,k}w_{2,l}-(l-1)w_{1,l}w_{2,k})r$, where $r\in R$ and 
$w_{t,i}^{(s)}=0$ for $s>0$, $t=1,2$, and $w_{t,i}=0$ for $i\ne k \vee  i\ne l$. 
It is not possible. In turn, 
\ba
c:K[1/(w_{1,k}w_{2,l}-w_{1,l}w_{2,k})][\{w_{1,i},w_{2,i},\dot{w}_{1,i},\dot{w}_{2,i},\ldots \}_{i=0}^{N}]\mapsto  \\  \nonumber
K[1/(w_{1,k}w_{2,l}-w_{1,l}w_{2,k})][\{w_{1,i},w_{2,i},\dot{w}_{1,i},\dot{w}_{2,i},\ldots \}_{i=0}^{N}]/I_{\max}=:\tilde{R}
\ea 
realizes a differential fields extension $K\subset \tilde{K}:=Q(R/I_{\max})$, where $Q(\ )$ is the field of fractions of a domain. 
$\overline{\tilde{K}}$ does not possess any PV extension 
via (\ref{pq}). Really, the finite set $\{c(w_{1,i}),c(w_{2,i})\}_{i=0}^N$ 
defines $w_1^c(y)/w_2^c(y)\in \tilde{K}(y) - \tilde{K}$ 
such that
\be \label{cc}
(w_1^c(y)/w_2^c(y))^{\bullet}|_{\dot{y}=P(y)/Q(y)}=0.
\ee
Let $K\subset K'(\subset \overline{\tilde{K}})$ be a maximal PV extension 
of $K$ via (\ref{pq}) in $\overline{\tilde{K}}$. Moreover, I assume 
that $\overline{K'}$ has a PV extension via (\ref{pq}) to 
$(\overline{K'}(y),\dot{y}=P(y)/Q(y))$. Therefore, 
$y\notin \overline{\tilde{K}}$, but $y$ fulfils (\ref{cc}) and its integral, 
so is algebraic over $\tilde{K}$. The contradiction ends the proof.

\end{proof}
One needs to remark that
the above $\tilde{K}$ is not PV extension of $K$, so $K$ does not possess 
necessary an integral of the type $\alpha (y-y_1)/(y-y_2)$.

Additionally, it is true that:
\begin{co}  \label{min}
Let $K_{\min}:= C\langle P_i,Q_j, \dot{P_i},\dot{Q_j},\ldots
\rangle
\subset K $, $P_i, Q_j$ are coefficients of
$P, Q \in K[y]$. Then $K_{\min}^{PQ}= K_{\min}\langle
y_{sol}\rangle $, where $\{y_{sol} \}\in K^{PQ}$ is the set of
solutions of (\ref{pq}) and $L^{PQ}$ is the greatest PV extension
via the equation of any differential field $L\supset K_{\min}$.
\end{co}
An inview into the structure of the arising set brings an example.
\begin{example}    \label{2}
Let the equation has the form
\be \label{riccati}
\dot{y}=y^2+by+a
\ee
(for $K=\bbold{C}$ this is Riccati equation). Then by the
transformation $y=-\frac{\dot{z}}{z}$ it becomes equivalent to
\be \label{linear}
\ddot{z}-b\dot{z}+az=0.
\ee
Let $z_1, z_2$ be linear independent solutions of the equation from
the PV extension $K^{PQ}$ via (\ref{linear}). Then all solutions of
(\ref{riccati}) are defined by the following algebraic set from
$K^{PQ}\times C$:
\be \label{alsets}
( (z_2y-\dot{z_2}) c+(z_1y-\dot{z_1}) )
\left(y-\frac{\dot{z_2} }{z_2} \right)=0.
\ee
Additionally, the greatest PV extension of $K$,
$\widetilde{K}\supset K $, via (\ref{riccati}) may be extended to
$K^{PQ}$ via $\dot{z}+y_s z=0$, where $y_s\in
\widetilde{K}$ is a solution of (\ref{riccati}). Therefore, the
questions about quadratures are equivalent for the two equations.
\end{example}
It is appeared that the algebraicity is a general feature of
the greatest sets of solutions.
\begin{pr}         \label{algebra}
Let $L\equiv K^{PQ} \supset K$ be the greatest PV extension of $K$ via
(\ref{pq}). Then the set of solutions is the
intersection of an algebraic set in $L \times
L$ and $\{(y,z)|y\in L \wedge z\in C \}$
subtracted a finite set.
\end{pr}
\begin{proof}
From lemma \ref{kac} one concludes that the extension of
differential fields $\overline{L}(y)\supset K$ contains NIC,
where $\dot{y}=\frac{P(y)}{Q(y)}$. Let $\vec{k}\in ({\cal
F}_0(\{y_{sol}\},\bbold{Z} ),+)$, $\vec{k}\ne 0$, be an integer
value function of finite support defined on the set of all
solutions.
The following
equation for $(y,c)\in \overline{L}\langle \alpha_{\vec{k}}\rangle \times C$ contain  the solutions set of (\ref{pq}):
\be \label{set}
\left( \alpha_{ \vec{k} }  \prod_{l,k_l>0}
(y-y_l)^{k_l}-c\prod_{l,k_l<0} (y-y_l)^{-k_l} \right)
\left( \prod_{l,k_l<0}
(y-y_l)^{-k_l} \right) = 0.
\ee
Moreover, the characteristic of $K$ is zero, so $\bbold{Q}\subset C$ and there are
infinite number of solutions of (\ref{pq}), roots of (\ref{nicform}). Therefore,
\be   \label{alphain}
\alpha_{\vec{k}}=c_0\prod_l(y_0-y_{l})^{-k_l}\in L
\ee
for a constant $c_0$ and a solution $y_0$. We may restrict to $L\times C$.

(\ref{set}) is equivalent to
\be     \label{setd}
\prod_l (y-y_l) \left( \frac{\dot{\alpha_{ \vec{k} } } }{\alpha_{\vec{k}} }
+\sum_l k_l \frac{1}{y-y_l} ( \dot{y}-\dot{y}_l ) \right) =0.
\ee
Therefore, if $y$ belongs to solutions of
(\ref{set}) for a $c\in C$ then $y$ is a solution of (\ref{pq}) or
\be  \label{ghost}
\sum_l k_l \frac{1}{y-y_l}=0.
\ee
\end{proof}
The common part of all admitted sets (\ref{set}) leads as a result to a
set $\Sigma_{PQ}$ in $\{(y,z)\in L\times L
\}$ intersected with $\{(y,z)| \dot{z}=0 \}$ covering all solutions
and a finite number of {\it ghost solutions} $y_g$ which are not
solutions of the differential equation and fulfil all additionally
arising equations (\ref{ghost}). $L[y,z]$ is a Noether
ring. It implies that a finite set of equations (\ref{ghost})
determined by a set $\{ \vec{k}_f \}_{f\in 1}^F$ may be chosen. The
following equations is also valid for ghost solutions
\be    \label{gadd}
\frac{ \dot{ \alpha_{\vec{k} } } }{\alpha_{\vec{k} }
}=\sum_l\frac{k_l\dot{y_l} }{y-y_l},
\ee
which follows from (\ref{setd}). In the way (\ref{gadd}) and
(\ref{ghost}) for $\vec{k}=\vec{k}_1,\ldots ,\vec{k}_F$ contains
all ghost solutions.

The ghost solutions appearance is natural and immovable. Consider a formal differential
complex: $d:K[y]\longrightarrow \Lambda K[y] \longrightarrow \Lambda^2 K[y]$
with the differentiation defined by: $dp(y):=\dot{p}(y)dx+p_y(y)dy$ for $p(y)
\in K[y]$ and $d(a(y)dx+b(y)dy):= (-a_y(y)+\dot{b}(y))dx\wedge dy$.
The notation is such that: $\dot{p}(y)\equiv \sum \dot{p_i}y^i$ and
$p_y(y)\equiv \sum ip_iy^{i-1}$ for $p(y)\in K[y]$. We have that
$d^2=0$, $H^0=C$ and $H^2=0$, where $H^i$ are $C$-linear spaces of cohomology.
If $\forall (k\in K) \int k \in K$ then we come
to the classical Poincare lemma: $H^1=0$. In turn, the Eq. (\ref{pq}) may be identify
with $\omega=P(y)dx-Q(y)dy=0$. Then the necessary integrability condition is $P_y(y)+
\dot{Q}(y)=0$ i.e. $d\omega
= 0$. Ghosts appear if an integral factor is needed.

The equation for an integral factor in $K(y)$, $\dot{y}=P(y)/Q(y)$, has the form:
\be
\dot{f}=-\frac{P_y(y)+\dot{Q}(y)}{Q(y)}f.
\ee
Therefore, if it exists it is unique up to a constant of $K(y)$.

For general type of equations
\be \label{pp}
p(y,\dot{y})=0,
\ee
$p\in K[y,\dot{y}]$ and $p_{\dot{y} }\ne 0$ the situation is
similar iff the roots $z$ of $p(y,z)=0$ are in $K(y)$. Otherwise
one concludes the following statement.
\begin{pr}
Let $p(y,z)=0$ has at least one root in $\overline{K(y)}-K(y)$.
Then for any set $\{y_i \}$ the extension $( K(\{ y_i \}) ,\dot{\ }
)\supset (K,\dot{\ })$ defined with agreement with
$p(y_i,\dot{y_i})=0$ is a PV extension of $K$ via $p(y_i,
\dot{y_i})=0$.
\end{pr}
\begin{proof}
The extension is defined by taking for each $\dot{y_i}$ a
nonrational root of $p$ in $\overline{K(y_i)}-K(y_i)$. If new
constants in $\overline{K(\{y_i\})}$ exist then a new constant
$w\in K( \{ y_i \} ) $ must belong to $K(y_1,\ldots ,y_N)$ for
$N\in \bbold{N}$. Then $w_{y_i}\ne 0$ for an $i\in [1,N]$, so
$\dot{y_i}\in \overline{K(\ldots,y_{i-1},y_{i+1},
\ldots )}[y_i]$. The contradiction.
\end{proof}
A formal series
analog of the algebraic constant of lemma \ref{kac}  may be found.
\begin{pr}  \label{wsol}
For $P,Q\in K[[y]]$, $P(y)=p_0y^{\mu_p}+\ldots$ and $Q(y)=q_0
y^{\mu_q}+
\ldots$, $p_0\ne 0$, $q_0=1$,
there is a PV extension $(\tilde{K},\dot{\ })\supset (K,
\dot{\ }) $ such that the equation:
\be \label{w}
\dot{w}Q+w_yP=0
\ee
has a nonzero solution $w=w_0y^{\mu_0}+\ldots$, $w_0\ne 0$, in
$\tilde{K}[[y]][y^{-1}]$. Moreover
\begin{enumerate}
\item[A] for $\mu_q+1\le \mu_p$ $w_k$ belongs to a liouvillian extension of
the differential field $ K_{\min} \langle \{w_l \}_{l=0}^{k-1},
\{\dot{w_l} \}_{l=0}^{k-1},
\ldots \rangle $, $K_{\min}:=
C\langle p_i,q_j,\dot{p_i},\dot{q_i},\ldots \rangle $, for all $
k\ge 0 $
\item[B] for $\mu_q+1 > \mu_p$ $w_k$ belongs to
the differential field $K_{\min}\langle \{ w_l \}_{l=0}^{k-1},
\{\dot{w}_{l} \}_{l=0}^{k-1}, \ldots
\rangle $ for all $k>0$.
\end{enumerate}
\end{pr}
\begin{proof}
The used method is a direct substitution. At the beginning, let
$\mu_0\in \bbold{Z}-\{0 \}$ and $w=x^{\mu_0}(w_0+w_1y+\ldots)$,
$w_0\ne 0$. Then necessary $\mu_q +1 < \mu_p $ is kept and then
\ba \label{ww}
q_0 \dot{w_0}=0 \\ \nonumber
\vdots  \\  \nonumber
q_0 \dot{ w_{\delta +k} } + \ldots + q_{\delta +k} \dot{w_0}
      +(\mu_0+k) p_0 w_k + \ldots +\mu_0 p_k w_0=0
          \\ \nonumber
\vdots ,
\ea
where $\delta := \mu_p-\mu_q -1$, or $\mu_q+1 = \mu_p$ and
\ba \label{ww=}
q_0 \dot{w_0} + \mu_0 p_0 w_0 = 0 \\ \nonumber
\vdots  \\  \nonumber
q_0 \dot{w_{ k } } +\ldots + q_{k}\dot{w_0}
      +(\mu_0+k)p_0 w_k + \ldots +\mu_0 p_k w_0=0
          \\ \nonumber
\vdots.
\ea
Therefore, we are in case A. Now, let $\mu_0=0$. If $\gamma :=
\mu_q+1- \mu_p > 0$ then
\ba \label{wwB}
q_0 \dot{w_0} + \gamma p_0 w_{\gamma } = 0 \\ \nonumber
\vdots  \\  \nonumber
q_0 \dot{w_{ k } } + \ldots +q_{k} \dot{w_0}
      +(\mu_0+k+\gamma)p_0 w_{k+\gamma} + \ldots +(\mu_0+\gamma) p_k w_{0+
\gamma}=0
          \\ \nonumber
\vdots
\ea
and $w_l=0$ for $l\in (0,\delta )$. It is case B. The rest
conditions for $\mu_p,\mu_q$ are of type A and are verified similarly.
\end{proof}
I return to the rational equation (\ref{pq}). Taking into the place
$w=w_0y^{\mu_0}+\ldots $ a solution of (\ref{w}) in the form of
(\ref{nicform}), (\ref{w}) is then (\ref{nic}), one obtains:
\begin{pr}        \label{liouv}
The greatest PV extension of $K$, $(K,\dot{\ })\subset
(K^{PQ},\dot{\ })$, via (\ref{pq}) is in an extended liouvillian
extension
(i.e. via
a finite sequence of extensions via an affine differential equation and an extension contained in
algebraic closure)
iff $K\subset K\langle w_0, \dot{w_0},\ldots
\rangle$ is in an extended liouvillian one, where $w=w_0y^{\mu_0}+\ldots
\notin K$ is a rational solution of (\ref{w}).
\end{pr}
\begin{proof}
Assume that
\be
w=\alpha \frac{a(y)}{b(y)}:=
\alpha y^{\mu_0} \frac{1+a_1y+\ldots +a_Ny^N}{1+\ldots +b_My^M}
\in L(y)-L
\ee
is a solution of (\ref{w}) and $K\subset K\langle w_0,\dot{w_0},
\ldots \rangle \subset L$ is extended liouvillian
and $a(y)$, $b(y)$ have no
common divisors, $\alpha \ne 0$. Its existence follows from lemma
\ref{kac} and theorem \ref{max}. Then a finite set of coefficients
of $w=w_0y^{\mu_0}+
\ldots $  define $\alpha, a_i, b_j$ by  algebraic
equations. Using proposition \ref{wsol} one states that they are in
an extended liouvillian extension $L'\supset K\langle w_0, \dot{w}_0
,\ldots
\rangle $. Moreover, $\{y|\exists_{c\in C}\ a(y) b(y) = c b(y)^2 \}$
consists from all solutions of (\ref{pq}) and a finite set of ghost
solutions, proposition \ref{algebra}.
Therefore, $K^{PQ}$ lies in an extended liouvillian extension of
$K$.

The inverse implication is obvious.
\end{proof}
Finally, from propositions \ref{wsol}, \ref{algebra} and the proof of 
proposition \ref{liouv}  we reach the following fact:
\begin{tm}  \label{end}
If $\mu_q+1\le \mu_p$ then $K\subset K^{PQ}$ is an extended
liouvillian extension.
\end{tm}


\bibliographystyle{amsplain}

\bibliography{nonlin}

\providecommand{\bysame}{\leavevmode\hbox to3em{\hrulefill}\thinspace}
\begin{thebibliography}{1}

\bibitem{Ar97}
J.~M. Aroca, \emph{Singularities of differential equations},  (1997),
  manuscript version.

\bibitem{Ar00}
\bysame, \emph{Puiseux solutions of singular differential equations},
  Resolution of singularities (Basel) (H.~Hauser et~al., ed.), Prog. Math.,
  vol. 181, Birkhauser, 2000, A research textbook in tribute to Oscar Zariski.
  Based on the courses given at the working week in Obergurgl, Austria,
  September 7-14, 1997, pp.~129--145.

\bibitem{BT98}
J.~W. Bruce and F.~Tari, \emph{Implicit differential equations from the
  singularity theory viewpoint}, Banach Center Publ., vol.~33, pp.~23--38,
  Polish Acad. Sci., Warsaw, 1996.

\bibitem{Ma94}
A.~G. Magid, \emph{{Lectures on Galois differential theory}}, Providence, AMS,
  1994.

\bibitem{Ma96}
T.~Maszczyk, \emph{Newlander-{N}irenberg type theorem for analytic algebras},
  Complex analysis and geometry (New York) (Vincenzo~Ancona et~al., ed.), Lect.
  Notes Pure Appl. Math., vol. 173, Marcel Dekker, 1996, Proceedings of the
  conference held in Trento, Italy, June 5-9, 1995, pp.~319--348.

\bibitem{Ma97}
\bysame, \emph{{Generalization of Newlander-Nirenberg theorem}}, PhD thesis
  (1997), (Polish).

\bibitem{Zo99}
H.~\.Zo\l{}\c{a}dek, \emph{{Monodromy Group}}, 1999, manuscripte version.

\end{thebibliography}

\end{document}